\input amstex
\documentstyle{amsppt}
\magnification=\magstep1

\nologo \NoBlackBoxes

\vsize 21.8 true cm \hsize 15.2 true cm \voffset 1 true cm
\lineskip=3pt plus 1pt minus 1pt

\leftheadtext{C. Baak, M. S. Moslehian} \rightheadtext{$\theta$-derivations on $JB^*$-triples}

\topmatter

\title On the Stability of $\theta$-Derivations on $JB^*$-Triples
\endtitle

\author Choonkil Baak$^1$, and Mohammad Sal Moslehian$^2$
\endauthor

\address
$^{1}$Department of Mathematics, Chungnam National University, Daejeon
305--764, South Korea; $^2$Department of Mathematics, Ferdowsi University, P.O. Box 1159,
Mashhad 91775, Iran
\endaddress

\email $^1$cgpark\@cnu.ac.kr; $^2$moslehian\@ferdowsi.um.ac.ir
\endemail

\abstract We introduce the concept of $\theta$-derivations
on $JB^*$-triples, and prove the Hyers--Ulam--Rassias
stability of $\theta$-derivations on $JB^*$-triples.
\endabstract

\thanks The first author was supported by Korea Research Foundation Grant
KRF-2005-070-C00009.
\endthanks

\keywords  Hyers--Ulam--Rassais stability, $\theta$-derivation,
$JB^*$-triple
\endkeywords

\subjclass 39B52, 39B82, 47B48, 17Cxx
\endsubjclass
\endtopmatter

\document

\baselineskip=0.75 true cm

\head 1.  Introduction \endhead

The original motivation to introduce the class of nonassociative
algebras known as Jordan algebras came from quantum mechanics (see
\cite{28}). Let $\Cal H$ be a complex Hilbert space, regarded as
the ``state space'' of a quantum mechanical system. Let $\Cal
L(\Cal H)$ be the real vector space of all bounded self-adjoint
linear operators on $\Cal H$, interpreted as the (bounded) {\it
observables} of the system. In 1932, Jordan observed that $\Cal
L(\Cal H)$ is a (nonassociative) algebra via the {\it
anticommutator product} $x \circ y : =\frac{xy+yx}{2}$. This is a
typical example of a (special) Jordan algebra. A commutative
algebra $\Cal B$ with product $x\circ y$ (not necessarily given by an
anticommutator) is called a {\it Jordan algebra} if $x^2\circ
(x\circ y) = x \circ (x^2\circ y)$ holds for all $x, y \in \Cal B$.

A complex Jordan algebra $\Cal B$ with a product $x \circ y$, and
a conjugate-linear algebra involution $x \mapsto x^*$ is called a
$JB^*$-algebra if $\Cal B$ carries a Banach space norm $\| \cdot
\| $ satisfying $\|x\| = \|x^*\|, \|x\circ y\| \le \|x\| \cdot
\|y\|$ and $\|\{xx^*x\}|| = \| x\|^3$ for all $x, y \in \Cal B$.
Here $\{xyz\}: = (x \circ y) \circ z + (y \circ z) \circ x - (x
\circ z) \circ y$ denotes the {\it Jordan triple product} of $x,
y, z \in \Cal B$ (see \cite{21, 22}).

The Jordan triple product of a $JB^*$-algebra leads us to a more
general algebraic structure, the so-called $JB^*$-triple, which
turns out to be appropriate for most applications to analysis. By
a (complex) $JB^*$-triple we mean a complex Banach space $\Cal J$
with a continuous triple product $$\{\cdot, \cdot, \cdot \} : \Cal J \times \Cal
J \times \Cal J \to \Cal J$$ which is linear in the outer
variables and conjugate linear in the middle variable, and has the
following properties:

(i) (commutativity) $\{x, y, z\} = \{z, y, x\}$;

(ii) (Jordan identity) $$L(a, b)\{x, y, z\} = \{L(a, b)x, y, z\} -
\{x, L(b, a)y, z\} + \{x, y, L(a, b)z\}$$ for all $a, b, x, y, z,
\in \Cal J$ in which $L(a, b)x : = \{a, b, x\}$;

(iii) For all $a \in \Cal J$ the operator $L(a, a)$ is hermitian,
i.e. $\|e^{itL(a, a)}\| = 1$, and has positive
spectrum in the Banach algebra $B(\Cal J)$;

(iv) $\|\{x, x, x\}\| = \|x\|^3$ for all $x \in \Cal J$.

The class of $JB^*$-triples contains all $C^*$-algebras via $\{x,
y, z\} = \frac{xy^*z + zy^*x}{2}$. Every $JB^*$-algebra is a
$JB^*$-triple under the triple product $$\{x, y, z\} : = (x \circ
y^*) \circ z + (y^* \circ z) \circ x - (x \circ z) \circ y^* .$$
Conversely, every $JB^*$-triple $\Cal J$ with an element $e$
satisfying $\{e, e, z\} = z$ for all $z \in \Cal J$, is a unital
$JB^*$-algebra equipped with the product $x \circ y : = \{x, e,
y\}$ and the involution $x^* : = \{e, x, e\}$; cf. \cite{9, 20,
26}.

The stability problem of functional equations originated from a
question of S.M. Ulam \cite{27} concerning the stability of group
homomorphisms: Let $(G_{1}, \ast)$ be a group and let
$(G_{2}, \diamond, d)$ be a metric group with the metric
$d(\cdot, \cdot)$. Given $\epsilon > 0$, does there exist a
$\delta (\epsilon)> 0$ such that if a mapping $h : G_1 \to G_2$
satisfies the inequality $d(h(x \ast y) , h(x) \diamond h(y)) < \delta$
for all $x,y \in G_1$, then there is a homomorphism $H: G_1 \to G_2 $ with
$d(h(x) , H(x)) < \epsilon $ for all $x \in G_1$?

If the answer is affirmative, we would say that the equation of
homomorphism $H(x \ast y)=H(x)\diamond H(y)$ is stable. The
concept of stability for a functional equation arises when we
replace the functional equation by an inequality which acts as a
perturbation of the equation. Thus the stability question of
functional equations is that how do the solutions of the
inequality differ from those of the given functional equation?

D.H. Hyers \cite{10} gave a first affirmative answer to the
question of Ulam in the context of Banach spaces with: Let $E_1$
and $E_2$ be Banach spaces.
Assume that $f : E_1 \to E_2$ satisfies $\|f(x+y)-f(x)-f(y)\| \le \epsilon$
for all $x, y \in E_1$ and some $\epsilon \ge 0.$ Then there
exists a unique additive mapping $T: E_1 \to E_2$ such that
$\|f(x)-T(x)\| \le \varepsilon$ for all $x \in E_1.$

Now assume that $E _1$ and $E_2$ are real normed spaces with
$E_2$ complete, $f:E_1\to E_2$ is a mapping such that
for each fixed $x\in E_1$, the mapping $t\mapsto f(tx)$ is
continuous on $\Bbb R$, and let there exist $\varepsilon\ge 0$ and
$p\neq 1$ such that
$$\|f(x+y)-f(x)-f(y)\|\le \varepsilon(\|x\|\sp p+\|y\|\sp p)$$ for all $x, y \in E_1$.

It was shown by Th. M. Rassias \cite{23} for $p\in [0, 1)$ (and
indeed $p<1$) and by Z. Gajda \cite{7} following the same approach
as in \cite{23} for $p>1$ that there exists a unique linear map
$T: E_1 \to E_2$ such that
$$\|f(x)-T(x)\|\le\frac{2\epsilon}{|2^p-2|}\|x\|\sp p$$ for all
$x \in E_1$. It is shown that there is no
analogue of Th.M. Rassias result for $p=1$ (see \cite{7, 25})

The inequality $\|f(x+y)-f(x)-f(y)\|\le \varepsilon(\|x\|\sp
p+\|y\|\sp p)$ has provided a lot of influence in the development
of what is now known as {\it Hyers--Ulam--Rassias stability} of
functional equations; cf. \cite{5, 6, 11, 13, 24}.

In 1992, G\u avruta \cite{8} proved the following.

\proclaim{Theorem 1.1}  Let $G$ be an abelian group and $X$ be a
Banach space. Denote by $\varphi : G \times G \to [0, \infty)$ a
function such that
$$\widetilde{\varphi}(x, y) = \sum_{j=0}^{\infty} \frac{1}{2^j}
\varphi(2^jx, 2^j y) < \infty$$
for all $x, y \in G$. Suppose that $f : G \to X$ is a mapping satisfying $$\|f(x+y) - f(x)-f(y) \|
\le \varphi(x, y)$$ for all $x, y \in G$. Then there exists a
unique additive mapping $T : G \to X$ such that
$$\|f(x)-T(x)\| \le \frac{1}{2}\widetilde{\varphi}(x, x)$$ for
all $x \in G$. \endproclaim

It is easy to see that Theorem 1.1 is still valid if
$$\widetilde{\varphi}(x, y) = \sum_{j=1}^{\infty} 2^{-j}
\varphi(2^{-j}x, 2^{-j}y) < \infty $$ (see also \cite{11}).

Since then the topic of approximate mappings or the stability of
functional equations was studied by several mathematicians;
\cite{2, 3, 15} and references therein. In particular, Jun and Lee
proved the following theorem; cf. \cite{12, Theorems 1 \& 6}.

\proclaim{Theorem 1.2} Denote by $\varphi : X  \times X \to [0,
\infty)$ a function such that
$$ \widetilde{\varphi}(x, y)=\sum_{j=0}^{\infty} \frac{1}{
3^j} \varphi(3^jx, 3^jy) < \infty \; \; \; \big(resp. \; \; \;
\widetilde{\varphi}(x, y)=\sum_{j=0}^{\infty} 3^j \varphi(3^{-j}x, 3^{-j}y) < \infty \big)$$
for all $x, y \in X $.
Suppose that $f : X \to Y$ is a mapping with $f(0) = 0$
satisfying $$\|2f(\frac{x+y}{2})-f(x)-f(y)\| \le \varphi(x, y)$$
for all $x, y \in X $. Then there exists a unique additive
mapping $T : X \to Y$ such that
$$\|f(x) - T(x)\| \le \frac{1}{3}\big(\widetilde{\varphi}(x, -x) +
\widetilde{\varphi}(-x, 3x)\big)$$
$$\big(resp. \; \; \; \|f(x) - T(x)\| \le \widetilde{\varphi}(\frac{x}{3}, \frac{-x}{3}) +
\widetilde{\varphi}(\frac{-x}{3}, x), \big)$$
for all $x\in X$.
\endproclaim

There are several various generalizations of the notion of
derivation. It seems that they are first appeared in the framework
of pure algebra (see \cite{1}). Recently they have been treated in
the Banach algebra theory (see \cite{14}). In addition, the
stability of these derivations is extensively studied by the
present authors and others; see \cite{4, 16, 18, 19} and
references therein.

In this paper, using some ideas from \cite{21}, we introduce the
notion of $\theta$-derivations on $JB^*$-algebras as a
generalization of derivations on $JB^*$-triples \cite{9} and prove
the Hyers--Ulam--Rassais stability of $\theta$-derivations on
$JB^*$-triples. Our result may be considered as a generalization
of those of \cite{20}.

\head 2.  Stability of $\theta$-derivations
\endhead

Throughout this section, let $\Cal J$ be a complex $JB^*$-triple with
norm $\| \cdot \|$.

\definition{Definition 2.1}
Let $\theta : \Cal J \to \Cal J$ be a $\Bbb C$-linear mapping. A $\Bbb C$-linear mapping $D :
\Cal J \to \Cal J$ is called a {\it $\theta$-derivation} on $\Cal J$ if
$$D(\{xyz\}) = \{D(x) \theta(y) \theta(z)\} + \{\theta(x)
D(y) \theta(z)\} + \{\theta(x) \theta(y) D(z)\}$$
for all $x, y, z \in \Cal J$.
\enddefinition

In particular, $D : = \frac{1}{3} \theta$ gives rise a
$JB^*$-homomorphism on $\Cal J$. Hence our results can be regarded
as an extension of those of \cite{20}. Note that if $D$ is a
derivation on a $JB^*$-algebra then every derivation $D$ can be
represented as $D_1 + i D_2$ where $D_1$ and $D_2$ are
$*$-preserving derivations.

\proclaim{Theorem 2.1}
 Let $f, h : \Cal J \to \Cal J$ be mappings with $f(0) =h(0)= 0$
 for which there exists a function $\varphi
: \Cal J^3 \to [0, \infty)$  such that
$$\align
\widetilde{\varphi}(x, y, z) :  = \sum_{j=0}^{\infty}
\frac{1}{2^{j}}
\varphi(2^j x, 2^j y, 2^j z)&   < \infty , \tag 2.1 \\
 \|f(\mu x + y) -  \mu f(x) - f(y) \|& \le \varphi(x, y, 0) ,  \tag 2.2 \\
 \|h(\mu x + y) - \mu h(x) - h(y)\| & \le \varphi(x, y, 0) , \tag 2.3 \\
 \|f(\{xyz\}) -  \{f(x)h(y)h(z)\} &  - \{h(x)f(y)h(z)\}
\\  - \{h(x) h(y) f(z)\}\| & \le \varphi(x, y, z) , \tag 2.4
\endalign$$ for all $x, y, z \in \Cal J$ and all $\mu \in S^{1} : = \{\lambda \in \Bbb C  \mid  | \lambda | = 1 \}$.
 Then there exist unique $\Bbb C$-linear
mappings $D, \theta : \Cal J \to \Cal J$ such that
$$\align \|f(x)- D(x)\| & \le \frac{1}{2} \widetilde{\varphi}(x, x, 0) , \tag 2.5 \\
\|h(x)- \theta(x)\| & \le \frac{1}{2} \widetilde{\varphi}(x, x, 0) \tag 2.6
\endalign $$
for all $x \in \Cal J$. Moreover, $D : \Cal J \to \Cal J$ is a $\theta$-derivation on $\Cal J$.
\endproclaim

\demo{Proof} Let $\mu = 1 \in S^1$ and $z=0$ in {\rm (2.2)} and
{\rm (2.3)}. It follows from Theorem 1.1 that there exist unique additive mappings $D, \theta : \Cal J \to \Cal J$
satisfying {\rm (2.5)} and {\rm (2.6)}. The additive mappings $D,
\theta : \Cal J \to \Cal J$ are given by
$$\align D(x) & = \lim_{l \to \infty} \frac{1}{2^l} f(2^l x) , \tag 2.7 \\
\theta(x) & = \lim_{l \to \infty} \frac{1}{2^l} h(2^l x)  \tag 2.8
\endalign$$
for all $x \in \Cal J$.

Let $\mu\in S^1$. Set $y = 0$ in {\rm (2.2)}. Then
$$\|f(\mu x) -\mu f(x)\| \le \varphi(x, 0, 0),$$
for all $x\in \Cal J$. So that
$$ 2^{-l}(f(\mu 2^l x) -\mu f(2^l x))\|\le 2^{-l} \varphi(2^l x, 0, 0),$$
for all $x\in\Cal J$. Since the right hand side tends to zero as $n\to\infty$, we have
$$ D(\mu x) = \lim_{l\to \infty}\frac{f(2^l \mu x)}{2^l}= \lim_{l\to
\infty}\frac{\mu f(2^lx)}{2^l} = \mu D(x)$$
for all $\mu\in\S^1$ and all $x\in\Cal J$.
Obviously, $D(0x)=0=0D(x)$.

Next, let $\lambda=\alpha_1+i\alpha_2\in \Bbb C$, where $\alpha_1,
\alpha_2\in \Bbb R$. Let
$\gamma_1=\alpha_1-\lfloor\alpha_1\rfloor,
\gamma_2=\alpha_2-\lfloor\alpha_2\rfloor$, in which $\lfloor
r\rfloor$ denotes the greatest integer less than or equal to the
number $r$. Then $0\le \gamma_i<1, ~(1\le i \le 2)$ and by using
Remark 2.2.2 of \cite{17} one can represent $\gamma_i$ as
$\gamma_i=\frac{\mu_{i,1} + \mu_{i,2}}{2}$ in which $\mu_{i,j}\in
S^1, ~~(1\le i,j \le 2)$. Since $D$ is additive we infer that
$$\align
D(\lambda x)&= D(\alpha_1x)+i D(\alpha_2x)\\
&=\lfloor\alpha_1\rfloor D(x)+D(\gamma_1x)+i\big( \lfloor\alpha_2\rfloor D(x)+D(\gamma_2x)\big )\\
&=\big( \lfloor\alpha_1\rfloor D(x)+\frac{1}{2}D(\mu_{1,1}x + \mu_{1,2}x)\big )+i \big ( \lfloor\alpha_2\rfloor D(x)+\frac{1}{2}D(\mu_{2,1}x + \mu_{2,2}x)\big )\\
&=\big(\lfloor\alpha_1\rfloor D(x)+\frac{1}{2}\mu_{1,1}D(x) +\frac{1}{2}\mu_{1,2}D(x)\big )\\
& +i(\lfloor\alpha_2\rfloor D(x)+\frac{1}{2}\mu_{2,1}D(x) + \frac{1}{2}\mu_{2,2}D(x)\big )\\
&= \alpha_1D(x)+i\alpha_2D(x)\\
&=\lambda D(x).
\endalign$$
for all $x\in \Cal J$. So that the additive mappings $D: \Cal J \to \Cal J$ is $\Bbb C$-linear.
A similar argument shows that $\theta$ is $\Bbb C$-linear.

It follows from {\rm (2.4)} that $$\align \frac{1}{2^{3l}}\| f(2^{3l}
\{xyz\}) - \{f(2^l x) h(2^ly)h(2^l z)\} & - \{h(2^l x) f(2^ly) h(2^l z)\}
 \\ - \{h(2^l x) h(2^l y) f(2^l z)\}\| \le \frac{1}{2^{3l}}
\varphi(2^l x, & 2^l y, 2^l z)  \le \frac{1}{2^l}
 \varphi(2^l x, 2^l y, 2^l z) , \endalign $$ which tends to
zero as $l \to \infty$ for all $x, y, z \in \Cal J$ by {\rm (2.1)}. By
{\rm (2.7)} and {\rm (2.8)},
$$ D(\{xyz\}) = \{D(x) \theta(y) \theta(z)\} + \{\theta(x) D(y) \theta(z)\}
+ \{\theta(x) \theta(y) D(z)\}$$ for all $x, y, z \in \Cal J$. So
the additive mapping $D : \Cal J \to \Cal J$ is a
$\theta$-derivation on $\Cal J$. \qed
\enddemo

\proclaim{Remark} {\rm It is easy to verify that the theorem is
true if
$$\widetilde{\varphi}(x, y) : = \sum_{j=1}^{\infty} 2^{-j} \varphi(2^{-j}x, 2^{-j}y) < \infty.$$}
\endproclaim

\proclaim{Corollary 2.2}  Let $f, h : \Cal J \to \Cal J$ be mappings with
$f(0) =h(0)= 0$ for which there exist constants $\epsilon \ge 0$
and $p \neq 1$ such that
$$\align \| f(\mu x + y)- \mu f(x) - f(y)\| & \le  \epsilon
(\|x\|^p + \|y\|^p),\\
\| h(\mu x +y)- \mu h(x)- h(y)\| & \le \epsilon (\|x\|^p + \|y\|^p),\\
\|f(\{xyz\}) - \{f(x)h(y)h(z)\} &  - \{h(x)f(y)h(z)\}\\
 - \{h(x) h(y) f(z)\}\|  & \le \epsilon (\|x\|^p + \|y\|^p + \|z\|^p) \endalign$$
for all  $x, y, z \in \Cal J$  and all $\mu \in S^{1}$.
Then there exist unique $\Bbb C$-linear mappings $D, \theta : \Cal J \to \Cal J$
such that
$$\align \|f(x)- D(x)\| & \le \frac{2\epsilon}{|2-2^p|} \|x\|^p ,  \\
\|h(x)- \theta(x)\| & \le \frac{2\epsilon}{|2-2^p|} \|x\|^p
\endalign $$
for all $x \in \Cal J$. Moreover, $D : \Cal J \to \Cal J$ is a $\theta$-derivation on $\Cal J$.
\endproclaim

\demo{Proof} Define $\varphi(x, y, z) = \epsilon (\|x\|^p +
\|y\|^p+\|z\|^p)$, and apply Theorem 2.1 and the remark following
the theorem. \qed
\enddemo

\proclaim{Theorem 2.3}
 Let $f, h : \Cal J \to \Cal J$ be mappings with
$f(0)=h(0) = 0$ for which there exists a function $\varphi : \Cal J^3
\to [0, \infty)$ satisfying {\rm (2.4)} such that
$$\align \widetilde{\varphi}(x, y, z) : = \sum_{j=0}^{\infty} \frac{1}{3^j}
\varphi(3^j x, 3^j y, 3^j z)& < \infty, \\
 \|2f(\frac{\mu x + y}{2}) - \mu f(x) - f(y)\| & \le \varphi(x, y, 0) ,  \tag 2.9
 \\ \|2h(\frac{\mu x +y}{2}) - \mu h(x) - h(y) \| & \le \varphi(x, y, 0) \tag 2.10
\endalign$$ for all $x, y, z \in \Cal J$ and all $\mu \in S^{1}$.
Then there exist unique $\Bbb C$-linear mappings $D, \theta : \Cal J \to \Cal J$
such that
$$\align \|f(x)- D(x)\| & \le \frac{1}{3}\big(\widetilde{\varphi}(x, -x, 0) +
\widetilde{\varphi}(-x, 3x, 0)\big) , \tag 2.11 \\
\|h(x)- \theta(x)\| & \le \frac{1}{3}\big(\widetilde{\varphi}(x, -x, 0)
+ \widetilde{\varphi}(-x, 3x, 0)\big) \tag 2.12
\endalign $$
for all $x \in \Cal J$. Moreover, $D : \Cal J \to \Cal J$ is a $\theta$-derivation on $\Cal J$.
\endproclaim

\demo{Proof} Let $z=0$ in {\rm (2.9)} and {\rm (2.10)}. It
follows from Theorem 1.2 that
there exist unique additive mappings $D, \theta : \Cal J \to \Cal J$
satisfying {\rm (2.11)} and {\rm (2.12)}. The additive mappings
$D, \theta : \Cal J \to \Cal J$ are given by
$$D(x) = \lim_{l \to \infty} \frac{1}{3^l} f(3^l x),$$
$$\theta(x) = \lim_{l \to \infty} \frac{1}{3^l} h(3^l x)$$
for all $x \in \Cal J$.

The rest of the proof is similar to the proof of Theorem 2.1.
\qed
\enddemo

\proclaim{Corollary 2.4}  Let $f, h : \Cal J \to \Cal J$ be mappings with
$f(0)=h(0) = 0$ for which there exist constants $\epsilon \ge 0$
and $p \in [0, 1)$ such that
$$\align \| 2f(\frac{\mu x + y }{2})- \mu f(x) - f(y)\| & \le  \epsilon
(\|x\|^p + \|y\|^p) ,  \\ \| 2h(\frac{\mu x + y}{2})- \mu h(x) - h(y)\| & \le \epsilon (\|x\|^p +
\|y\|^p),
 \\ \|f(\{xyz\}) -  \{f(x)h(y)h(z)\} &  - \{h(x)f(y)h(z)\}
\\  - \{h(x) h(y) f(z)\}\|  & \le
\epsilon (\|x\|^p + \|y\|^p + \|z\|^p)
\endalign$$ for all  $x, y, z \in \Cal J$ and all $\mu \in S^{1}$.
Then there exist unique $\Bbb C$-linear mappings $D, \theta : \Cal J \to \Cal J$
such that
$$\align \|f(x)- D(x)\| & \le \frac{3+3^p}{3-3^p} \epsilon \|x\|^p ,  \\
\|h(x)- \theta(x)\| & \le \frac{3+3^p}{3-3^p} \epsilon \|x\|^p
\endalign $$
for all $x \in \Cal J$. Moreover, $D : \Cal J \to \Cal J$ is a $\theta$-derivation on $\Cal J$.
\endproclaim

\demo{Proof} Define $\varphi(x, y, z) = \epsilon (\|x\|^p +
\|y\|^p + \|z\|^p)$, and apply Theorem 2.3. \qed
\enddemo

\proclaim{Theorem 2.5}
 Let $f, h : \Cal J \to \Cal J$ be mappings with
$f(0)=h(0) = 0$ for which there exists a function $\varphi : \Cal J^3
\to [0, \infty)$ satisfying {\rm (2.9)}, {\rm (2.10)} and {\rm
(2.4)}  such that
$$\sum_{j=0}^{\infty} 3^{3j}
\varphi(\frac{x}{3^j}, \frac{y}{3^j}, \frac{z}{3^j}) < \infty
\tag 2.16 $$ for all $x, y, z \in \Cal J$. Then there exist unique
$\Bbb C$-linear mappings $D, \theta : \Cal J \to \Cal J$ such that
$$\align \|f(x)- D(x)\| & \le \widetilde{\varphi}(\frac{x}{3}, - \frac{x}{3}, 0) +
\widetilde{\varphi}(-\frac{x}{3}, x, 0) , \tag 2.17 \\
\|h(x)- \theta(x)\| & \le \widetilde{\varphi}(\frac{x}{3}, - \frac{x}{3}, 0) +
\widetilde{\varphi}(-\frac{x}{3}, x, 0)  \tag 2.18
\endalign $$
for all $x \in \Cal J$, where $$ \widetilde{\varphi}(x, y, z) : =
\sum_{j=0}^{\infty} 3^j \varphi(\frac{x}{3^j}, \frac{y}{3^j},
\frac{z}{3^j})$$ for all $x, y, z \in \Cal J$. Moreover, $D : \Cal J \to \Cal J$
is a $\theta$-derivation on $\Cal J$.
\endproclaim

\demo{Proof} By Theorem 1.2, it follows from {\rm (2.16)}, {\rm (2.9)} and {\rm (2.10)} that
there exist unique additive mappings $D, \theta : \Cal J \to \Cal J$
satisfying {\rm (2.17)} and {\rm (2.18)}. The additive mappings
$D, \theta : \Cal J \to \Cal J$ are given by
$$\align D(x) & = \lim_{l \to \infty} 3^l f(\frac{x}{3^l}) , \tag 2.19 \\
\theta(x) & = \lim_{l \to \infty} 3^l h(\frac{x}{3^l}) \tag 2.20
\endalign
$$ for all $x \in \Cal J$.

By a similar method to the proof of Theorem 2.1, one can show
that $D, \theta : \Cal J\to \Cal J$ are $\Bbb C$-linear mappings.

It follows from {\rm (2.4)} that $$\align 3^{3l} \| f(\frac{\{xyz\}}{3^{3l}})
- \{f(\frac{x}{3^l}) h(\frac{y}{3^l})h(\frac{z}{3^l})\} & - \{h(\frac{x}{3^l})f(\frac{y}{3^l}) h(\frac{z}{3^l})\}
\\ -\{h(\frac{x}{3^l})h(\frac{y}{3^l}) f(\frac{z}{3^l})\}\| & \le 3^{3l}
\varphi(\frac{x}{3^l}, \frac{y}{3^l}, \frac{z}{3^l}) , \endalign
$$ which tends to zero as $l \to \infty$ for all $x, y, z \in \Cal J$ by
{\rm (2.16)}. By {\rm (2.19)} and {\rm (2.20)},
$$ D(\{xyz\}) = \{D(x) \theta(y) \theta(z)\} + \{\theta(x) D(y) \theta(z)\}
 + \{\theta(x) \theta(y) D(z)\}$$ for all $x, y, z \in \Cal J$.
So the additive mapping $D : \Cal J \to \Cal J$ is a
$\theta$-derivation on $\Cal J$. \qed
\enddemo

\proclaim{Corollary 2.6}  Let $f, h : \Cal J \to \Cal J$ be mappings with
$f(0) =h(0)= 0$ for which there exist constants $\epsilon \ge 0$
and $p \in (3, \infty)$ such that
$$\align \| 2f(\frac{\mu x +y}{2})-\mu f(x) -f(y)\| & \le  \epsilon
(\|x\|^p + \|y\|^p),  \\
\| 2h(\frac{\mu x + y}{2})- \mu h(x) - h(y)\| & \le \epsilon (\|x\|^p +
\|y\|^p),\\
\|f(\{xyz\}) -  \{f(x)h(y)h(z)\} &  - \{h(x)f(y)h(z)\}\\
- \{h(x) h(y) f(z)\}\| & \le
\epsilon (\|x\|^p + \|y\|^p + \|z\|^p)
\endalign$$
for all  $x, y, z \in \Cal J$ and all $\mu \in S^{1}$.
Then there exist unique $\Bbb C$-linear mappings $D, \theta : \Cal J \to \Cal J$
such that
$$\align \|f(x)- D(x)\| &  \le \frac{3^p+3}{3^p-3} \epsilon \|x\|^p ,  \\
\|h(x)- \theta(x)\| & \le \frac{3^p+3}{3^p-3} \epsilon \|x\|^p
\endalign $$
for all $x \in \Cal J$. Moreover, $D : \Cal J \to \Cal J$ is a $\theta$-derivation on $\Cal J$.
\endproclaim

\demo{Proof} Define $\varphi(x, y, z) = \epsilon (\|x\|^p +
\|y\|^p + \|z\|^p)$, and apply Theorem 2.5. \qed
\enddemo

\definition{Definition 2.2}
Let $\theta : \Cal J \to \Cal J$ be a $\Bbb C$-linear mapping. A $\Bbb C$-linear mapping $D :
\Cal J \to \Cal J$ is called a {\it Jordan $\theta$-derivation} on
$\Cal J$ if
$$D(\{xxx\}) = \{D(x) \theta(x)\theta(x)\} + \{\theta(x)
D(x) \theta(x)\} + \{\theta(x) \theta(x) D(x)\}$$ holds for
all $x \in \Cal J$.
\enddefinition

\remark{Problem 2.1} Is every Jordan $\theta$-derivation
 a $\theta$-derivation?
\endremark

\enddocument